\documentclass[14pt]{article}
\usepackage{amsfonts}
\usepackage{amsmath}
\usepackage{graphicx}

\newtheorem{theorem}{Theorem}

\newtheorem{proposition}[theorem]{Proposition}

\begin{document}

\title{\textbf{Synchronization and secure communication using some chaotic systems of
fractional differential equations}}
\author{ \textbf{O. Chi\c{s}$^*$, D. Opri\c{s}$^{**}$}}
\maketitle
\begin{center}
$^*, \, ^{**}$  Faculty of Mathematics and Informatics, West University of Timi\c{s}oara, Romania\\
E-mail: chisoana@yahoo.com, opris@math.uvt.ro
\end{center}

\textbf{Abstract}: Using Caputo fractional derivative of order
$\alpha,$ $\alpha\in (0,1),$ we consider some chaotic systems of
fractional differential equation. We will prove that they can be
synchronized and anti-synchronized using suitable nonlinear
control function. The synchronized or anti-synchronized error
system of fractional differential equations is used in secure communication.\\

\maketitle \textbf{MSC2000}: 65P20, 94A05, 11T71.\\


\section{Introduction}

Synchronization phenomenon has been studied intensively because of
its application in many fields, and one of it is secure
communication. There are many ways for synchronization, such as
feedback method, adaptive techniques, time delay feedback
approach, backstepping method, with nonlinear control \cite{Zha}.
We will prove here synchronization and anti-synchronization
between two chaotic systems of differential equations, by
considering a suitable nonlinear control function.

In the first section we will show synchronization between some
representative chaotic systems of fractional differential
equations, coupled fractional systems T and between system T and
R\"{o}ssler system. In Section 2 we will present
anti-synchronization of the same chaotic systems and we will
compare the two methods. In Section 3 we will apply
synchronization in secure communication. Numerical simulations are
done using Adams-Bashforth-Moulton algorithm \cite{Die}. In last
section some conclusions are presented.

We will briefly give the definition of fractional derivative, of
the following form
\begin{equation}\label{1}
D^\alpha_t x(t):=I^{m-\alpha}\Big(\frac{d}{dt}\Big)^m x(t), \,
\alpha>0,
\end{equation}
where $m=[\alpha],$
$\Big(\frac{d}{dt}\Big)^m=\frac{d}{dt}\circ...\circ \frac{d}{dt},$
$I^\beta$ is the $\beta-$order Riemann-Liouville integral operator
and it is expressed as
$$I^\beta_t x(t)=\frac{1}{\Gamma(\beta)}\int_0^t(t-s)^{\beta-1}x(s)ds, \, \beta>0,$$
where  $\Gamma$ is the Gamma function \cite{Die}.

Along this paper we will work with the following chaotic systems
of fractional diffe\-rential equations:
\begin{enumerate}
    \item T system of fractional differential equations \cite{Tig}
             \begin{equation}\label{T1}
             \left\{%
             \begin{array}{ll}
                D^{\alpha_1}_t x_1(t)=a_1(y_1(t)-x_1(t)),\\
                D^{\alpha_2}_t y_1(t)=(c_1-a_1)x_1(t)-a_1x_1(t)z_1(t),\\
                D^{\alpha_3}_t z_1(t)=x_1(t)y_1(t)-b_1z_1(t),\\
                \end{array}
                \right.
             \end{equation} where $\alpha_1, \alpha_2, \alpha_3 \in (0,1)$ and $a_1, \, b_1, \,
             c_1$ are the parameters of the system. It is known that the system has a chaotic
             behaviour for $a_1:=2.1,$ $b_1:=0.6$ and $c_1:=30;$
    \item Fractional version of R\"{o}ssler system with $\alpha_1, \alpha_2, \alpha_3 \in (0,1)$
    is
        \begin{equation}\label{R1}
         \left\{%
             \begin{array}{ll}
                D^{\alpha_1}_t x_1(t)=-y_1(t)-z_1(t),\\
                D^{\alpha_2}_t y_1(t)=x_1(t)+a_2y_1(t),\\
                D^{\alpha_3}_t z_1(t)=b_2+z_1(t)(x_1(t)-c_2)\\
                \end{array}
                \right.
             \end{equation} If  $a_2:=0.2,$ $b_2:=0.2$ and
             $c_2:=5.7,$ then the system is chaotic.
\end{enumerate}

\section{Synchronization between two chaotic fractional systems via nonlinear control}

The systems implied in synchronization are called drive (master)
and response (slave) systems. Let us consider the drive system
given in the form $D^{\alpha}_t x(t)=f(x(t))$ and the response
system $D^{\alpha}_t y(t)=g(x(t),y(t)),$ where $x(t)=(x_1(t),
x_2(t),x_3(t))\in \mathbb{R}^3$ and
$y(t)=(y_1(t),y_2(t),y_3(t))\in \mathbb{R}^3$ are the phase space
variables and $f, \, g$ the corresponding nonlinear functions. The
two systems will be synchronous if the trajectory of drive system
follows the same path as the response system, that means
\begin{equation} |x(t)-y(t)|\rightarrow c, \, t\rightarrow\infty.
\end{equation}
 If $c=0,$ then the synchronization is called
complete synchronization. Synchronization can be bone using drive
and response systems of the same type, or different types of
chaotic systems of fractional differential equations. We will
consider both these types of synchronization using the systems
(\ref{T1}) and (\ref{R1}) \cite{Li}, \cite{Seb}.

We will begin with the synchronization between two identical
fractional T systems (coupled) in the following manner. We
consider the drive system (\ref{T1}) and the response system with
a control $u(t)=[u_1(t),u_2(t),u_3(t)]^T$
\begin{equation}\label{T2}
             \left\{%
             \begin{array}{ll}
                D^{\alpha_1}_t x_2(t)=a_1(y_2(t)-x_2(t))+u_1(t),\\
                D^{\alpha_2}_t y_2(t)=(c_1-a_1)x_2(t)-a_1x_2(t)z_2(t)+u_2(t),\\
                D^{\alpha_3}_t z_2(t)=x_2(t)y_2(t)-b_1z_2(t)+u_3(t).\\
                \end{array}
                \right.
\end{equation}
The error system is given by subtracting (\ref{T2}) and
(\ref{T1}), and we get
\begin{equation}\label{Te1}
             \left\{%
             \begin{array}{ll}
                D^{\alpha_1}_t e_1(t)=a_1(e_2(t)-e_1(t))+u_1(t),\\
                D^{\alpha_2}_t e_2(t)=(c_1-a_1)e_1(t)-a_1(x_2(t)z_2(t)-x_1(t)z_1(t))+u_2(t),\\
                D^{\alpha_3}_t e_3(t)=-b_1e_3(t)+x_2(t)y_2(t)-x_1(t)y_1(t)+u_3(t),\\
                \end{array}
                \right.
\end{equation}where $e_1(t):=x_2(t)-x_1(t),$
$e_2(t):=y_2(t)-y_1(t),$ $e_3(t):=z_2(t)-z_1(t).$ We will re-write
the control $u$ by considering another control $v$ suitable chosen
that is a function of error states $e_1, e_2, e_3.$ We re-refine
$u$ as
\begin{equation}\label{c1}
\begin{array}{ll}
u_1(t):=v_1(t),\\
u_2(t):=a_1(x_2(t)z_2(t)-x_1(t)z_1(t))+v_2(t),\\
u_3(t):=x_1(t)y_1(t)-x_2(t)y_2(t)+v_3(t).
\end{array}
\end{equation}
Substituting (\ref{c1}) in (\ref{Te1}) we get
\begin{equation}\label{Te2}
             \left\{%
             \begin{array}{ll}
                D^{\alpha_1}_t e_1(t)=a_1(e_2(t)-e_1(t))+v_1(t),\\
                D^{\alpha_2}_t e_2(t)=(c_1-a_1)e_1(t)+v_2(t),\\
                D^{\alpha_3}_t e_3(t)=-b_1e_3(t)+v_3(t).\\
                \end{array}
                \right.
\end{equation}
We want to prove that the two considered systems are globally
synchronized, that means we have to choose the feedback control
$v$ such that the error to converge to 0, when
$t\rightarrow\infty.$ We choose $v$ such that
\begin{equation}\label{v}
[v_1(t),v_2(t),v_3(t)]^T=A\cdot[e_1(t),e_2(t),e_3(t)]^T,
\end{equation} where $A\in\mathcal{M}_{3\times3}.$ One choice of
$A$ is $$A:=\begin{bmatrix}
  0 & a_1 & 0\\
  -(c_1-a_1) & c_1 & 0\\
  0 & 0 & 2b_1 \\
\end{bmatrix}$$ and the feedback functions are
\begin{equation}\label{vT}
\begin{array}{ll}
v_1(t)=a_1e_2(t),\\
v_2(t)=-(c_1-a_1)e_1(t)+c_1e_2(t),\\
v_3(t)=2b_1e_3(t).
\end{array}
\end{equation}
The error system (\ref{Te2}) becomes
\begin{equation}\label{Te3}
             \left\{%
             \begin{array}{ll}
                D^{\alpha_1}_t e_1(t)=a_1e_1(t),\\
                D^{\alpha_2}_t e_2(t)=c_1e_2(t),\\
                D^{\alpha_3}_t e_3(t)=b_1e_3(t).\\
                \end{array}
                \right.
\end{equation}

\begin{proposition}\label{Pr1}
The systems of fractional differential equations \emph{(\ref{T1})}
and \emph{(\ref{T2})} will a\-pproach global asymptotical
synchronization for any initial conditions and with the feedback
control \emph{(\ref{vT})}.
\end{proposition}

\textbf{Proof:} We will prove asymptotic stability of the system
(\ref{Te3}) by using the Laplace transform \cite{Mu}. We take
Laplace transform in both sides of (\ref{Te3}), with
$E_i(s)=\mathcal{L}(e_i(t)), \, i=1,2,3,$ where
$\mathcal{L}(D^{\alpha_i}e_i(t))=s^{\alpha_i}E_i(s)-s^{\alpha_i-1}e_i(0),$
$i=1,2,3$ and we get
$$E_1(s)=\frac{s^{\alpha_1-1}e_1(0)}{s^{\alpha_1}-a_1}, \,
E_2(s)=\frac{s^{\alpha_2-1}e_2(0)}{s^{\alpha_2}-c_1}, \,
E_3(s)=\frac{s^{\alpha_3-1}e_3(0)}{s^{\alpha_3}-b_1}.$$ By
final-value theorem of the Laplace transformation \cite{Mu}, we
have
$$\mathop{\lim} \limits_{t\to\infty}e_1(t)=\mathop{\lim} \limits_{s\to 0}sE_1(s)=0, \,
\mathop{\lim} \limits_{t\to\infty}e_2(t)=\mathop{\lim}
\limits_{s\to 0}sE_2(s)=0, \, \mathop{\lim}
\limits_{t\to\infty}e_3(t)=\mathop{\lim} \limits_{s\to
0}sE_3(s)=0.$$ Therefore, systems (\ref{T1}) and (\ref{T2}) can
achieve asymptotic synchronization for any initial conditions and
with the control (\ref{vT}). \hfill $\Box$\\

For numerical simulation we use Adams-Bashforth-Moulton algorithm.
For $\alpha_1=0.9,$ $\alpha_2=0.5,$ $\alpha_3=0.6$ and initial
conditions $x_1(0)=0.01,$ $x_2(0)=0.01,$ $x_3(0)=0.01.$ Orbits of
the error system (\ref{Te3}) are represented in the above figures.

\begin{center}\begin{tabular}{ccc}
  \includegraphics[height=2in]{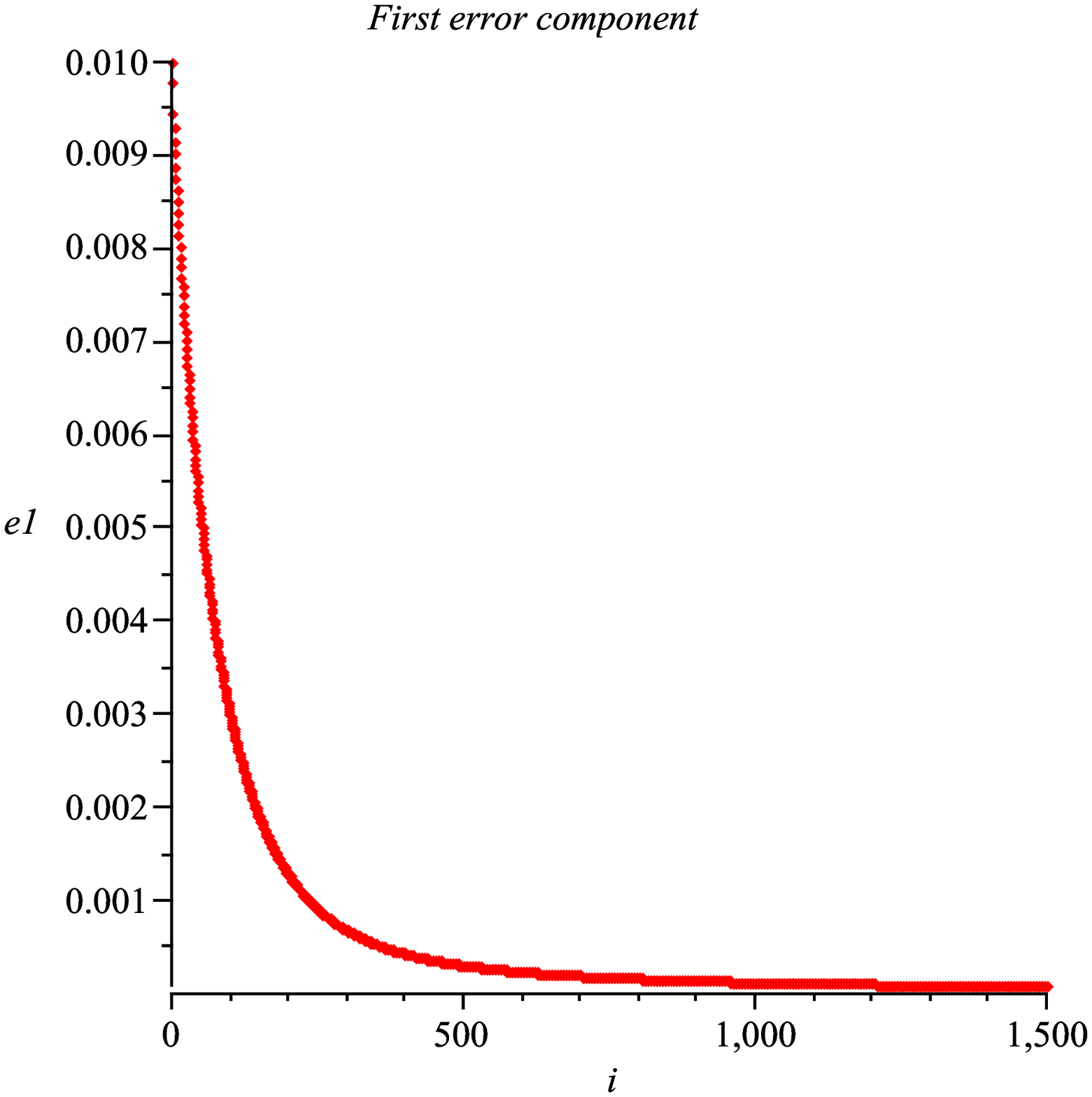}
  \includegraphics[height=2in]{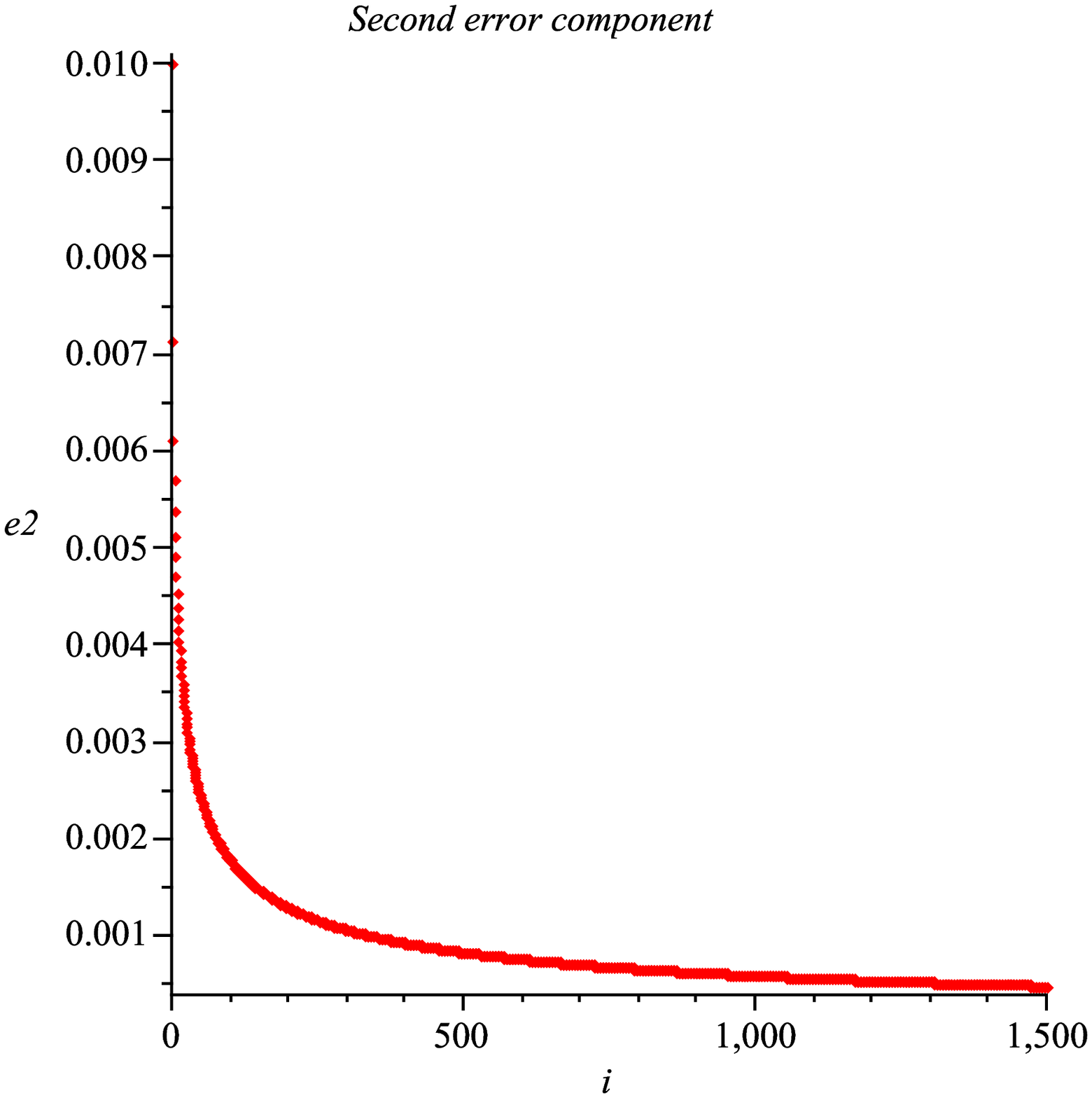}
  \includegraphics[height=2in]{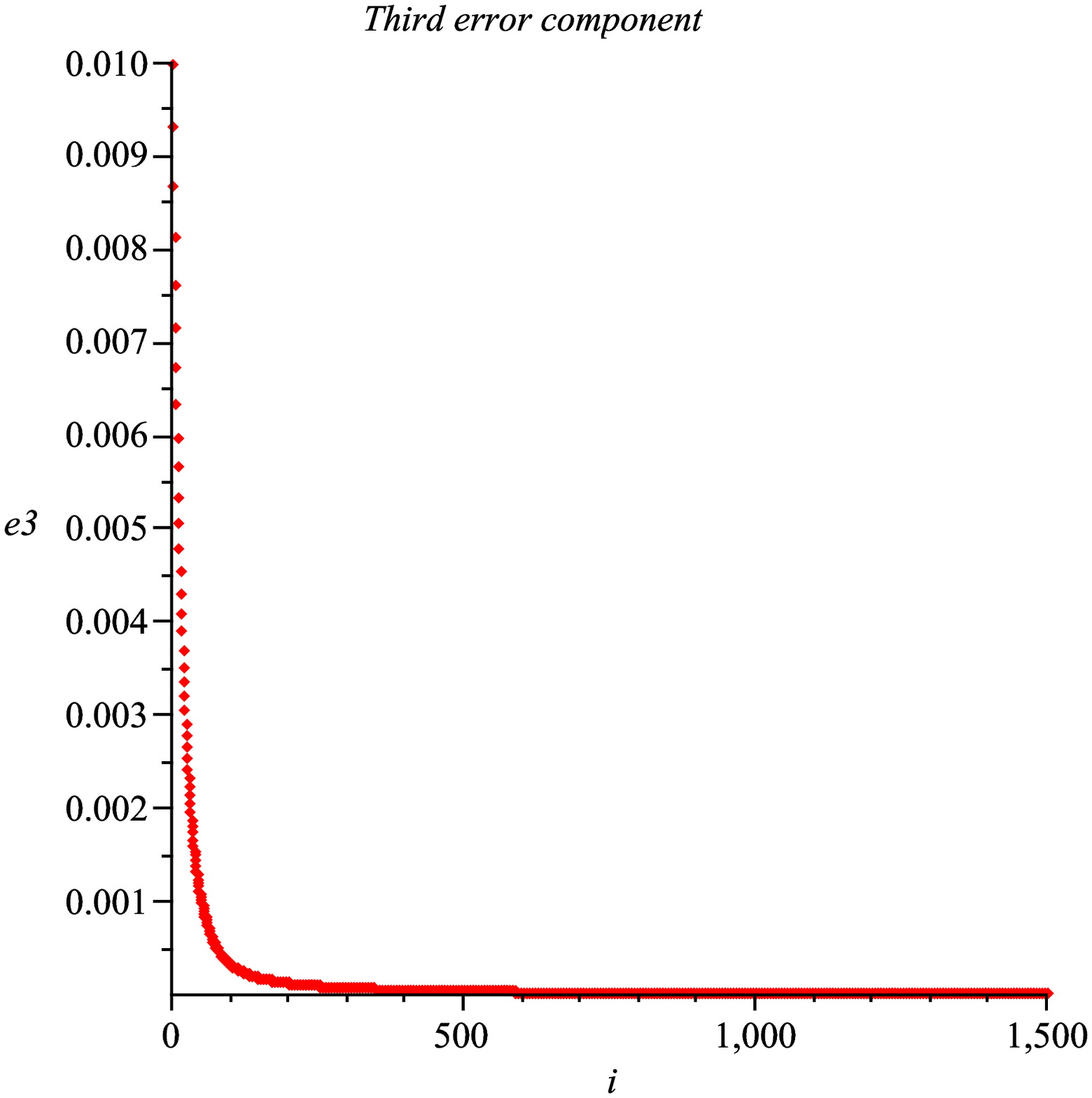}
\end{tabular}
\end{center}

We will consider synchronization between two different chaotic
systems of fractional differential equations, R\"{o}ssler system
as drive system (\ref{R1})
 and T system as response system
\begin{equation}\label{T11}
             \left\{%
             \begin{array}{ll}
                D^{\alpha_1}_t x_2(t)=a_1(y_2(t)-x_2(t))+u_1(t),\\
                D^{\alpha_2}_t y_2(t)=(c_1-a_1)x_2(t)-a_2x_2(t)z_2(t)+u_2(t),\\
                D^{\alpha_3}_t z_2(t)=x_2(t)y_2(t)-b_1z_2(t)+u_3(t).\\
                \end{array}
                \right.
             \end{equation}
The error system is obtained by subtracting (\ref{T11}) and
(\ref{R1}) and it has the form
\begin{equation}\label{Re1}
             \left\{%
             \begin{array}{ll}
                D^{\alpha_1}_t e_1(t)=-a_1e_1(t)+(a_1-1)e_2(t)-e_3(t)+a_1(y_1(t)-x_1(t))
                +y_2(t)+z_2(t)+u_1(t),\\
                D^{\alpha_2}_t
                e_2(t)=(c_1-a_1+1)e_1(t)+a_2e_2(t)+(c_1-a_1)x_1(t)-x_2(t)-
                a_2x_2(t)z_2(t)- \\ \quad \quad \quad \quad \quad \,  -a_2y_2(t)+u_2(t),\\
                D^{\alpha_3}_t
                e_3(t)=-(b_1+c_2)e_3(t)-b_2+x_2(t)y_2(t)-b_1z_1(t)-z_1(t)x_1(t)+
                c_2z_2(t)+u_3(t),\\
                \end{array}
                \right.
\end{equation}where $e_1(t):=x_2(t)-x_1(t),$
$e_2(t):=y_2(t)-y_1(t),$ $e_3(t):=z_2(t)-z_1(t).$ We re-write
control $u$ by considering the control $v$ such that
\begin{equation}
\begin{array}{ll}
u_1(t):=-a_1(y_1(t)-x_1(t))-y_2(t)-z_2(t)+v_1(t),\\
u_2(t):=-(c_1-a_1)x_1(t)+x_2(t)+a_2x_2(t)z_2(t)+a_2y_2(t)+v_2(t),\\
u_3(t):=b_2-x_2(t)y_2(t)+b_1z_1(t)+z_1(t)x_1(t)-c_2z_2(t)+v_3(t),
\end{array}
\end{equation}and the error system (\ref{Re1}) becomes
\begin{equation}\label{Re2}
             \left\{%
             \begin{array}{ll}
                D^{\alpha_1}_t e_1(t)=-a_1e_1(t)+(a_1-1)e_2(t)-e_3(t)+v_1(t),\\
                D^{\alpha_2}_t e_2(t)=(c_1-a_1+1)e_1(t)+a_2e_2(t)+v_2(t),\\
                D^{\alpha_3}_t e_3(t)=-(b_1+c_2)e_3(t)+v_3(t).\\
                \end{array}
                \right.
\end{equation}
As in the case above, we have to chose the feedback control $v$
such that the relation (\ref{v}) to be fulfilled. One choice of
matrix $A$ will be
$$A:=\begin{bmatrix}
  2a_1 & -(a_1-1) & 1 \\
  -(c_1-a_1+1) & 0 & 0\\
  0 & 0 & 2b_1+c_2 \\
\end{bmatrix}$$ and the feedback functions are
\begin{equation}\label{vR}
\begin{array}{ll}
v_1(t)=2a_1e_1(t)-(a_1-1)e_2+e_3(t),\\
v_2(t)=-(c_1-a_1+1)e_1(t),\\
v_3(t)=(2b_1+c_2)e_3(t).
\end{array}
\end{equation}
The error system (\ref{Re2}) becomes
\begin{equation}\label{Re3}
             \left\{%
             \begin{array}{ll}
                D^{\alpha_1}_t e_1(t)=a_1e_1(t),\\
                D^{\alpha_2}_t e_2(t)=a_2e_2(t),\\
                D^{\alpha_3}_t e_3(t)=b_1e_3(t).\\
                \end{array}
                \right.
\end{equation}

\begin{center}\begin{tabular}{ccc}
  \includegraphics[height=2in]{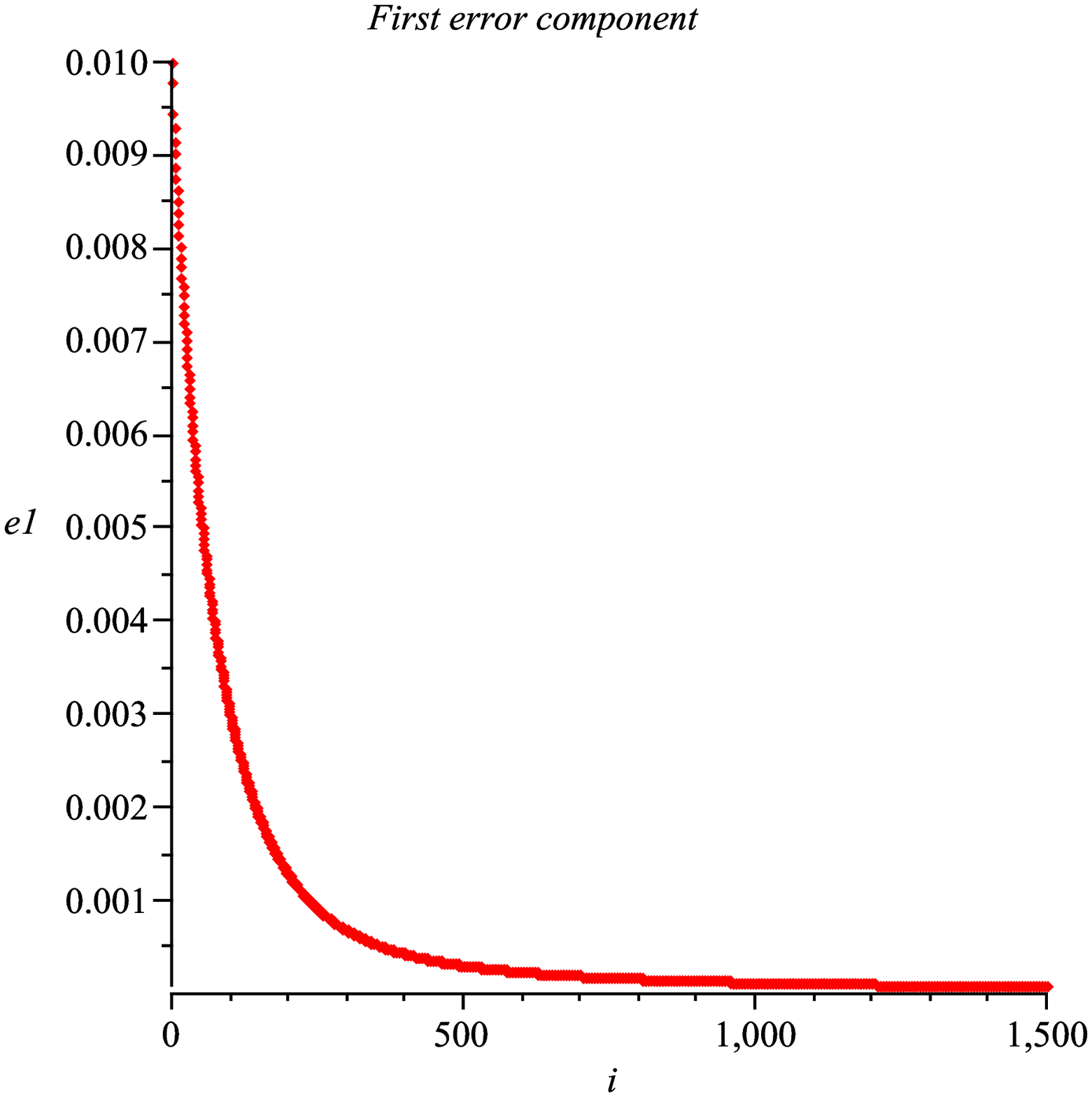}
  \includegraphics[height=2in]{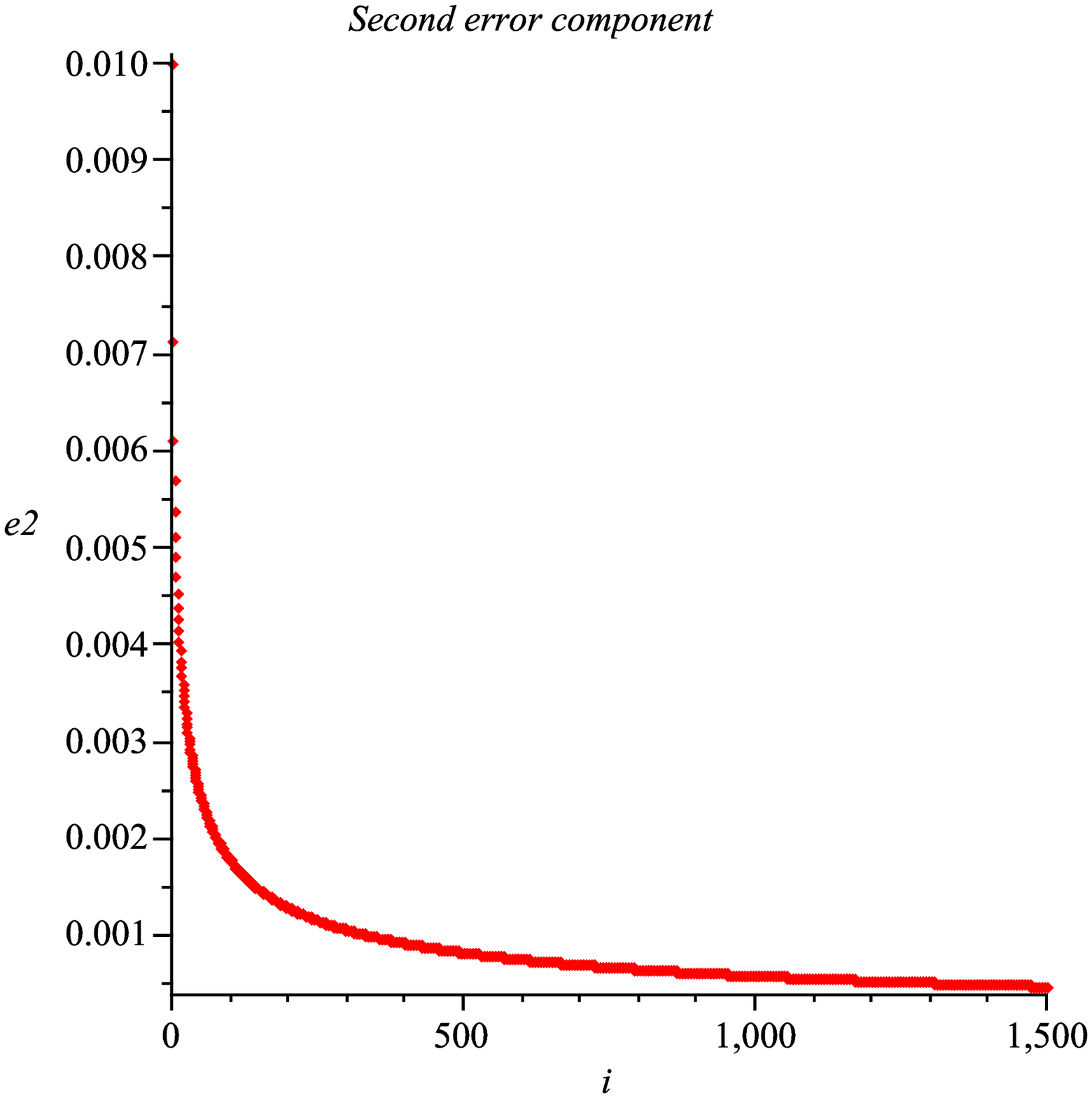}
  \includegraphics[height=2in]{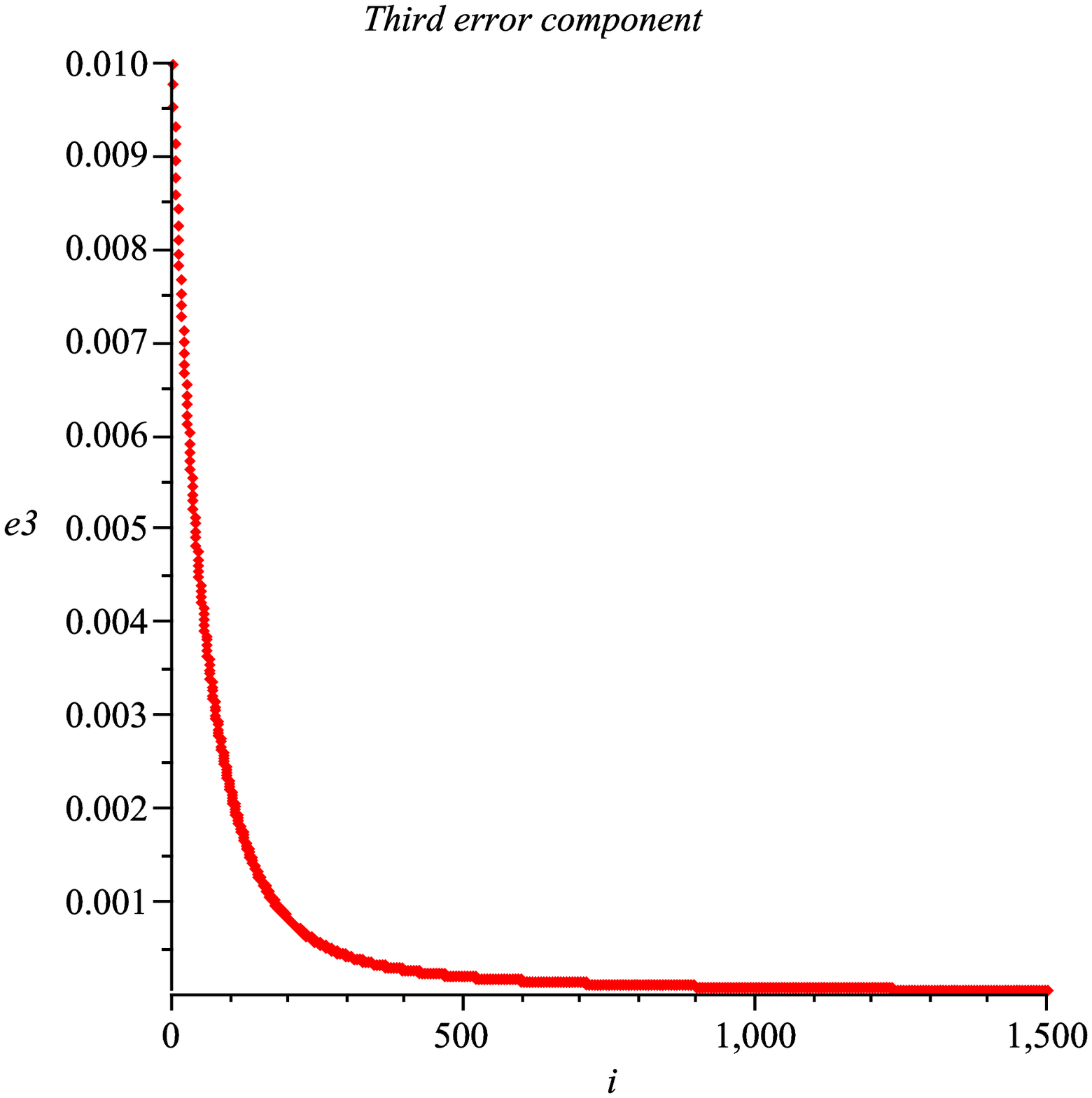}
\end{tabular}
\end{center}

\begin{proposition}
The systems of fractional differential equations \emph{(\ref{R1})}
and \emph{(\ref{T11})} will a\-pproach global asymptotical
synchronization for any initial conditions and with the feedback
control \emph{(\ref{vR})}.
\end{proposition}
\textbf{Proof:} The proof is similar to that of Proposition
\ref{Pr1}. \hfill $\Box$

The error system (\ref{Re3}) is represented in the above figures,
for $\alpha_1=0.9,$ $\alpha_2=0.5,$ $\alpha_3=0.6$ and initial
conditions $x_1(0)=0.01,$ $x_2(0)=0.01,$ $x_3(0)=0.01.$

\section{Anti-synchronization between two chaotic fractional systems via nonlinear control}

The drive and response  systems will be anti-synchronous if
$$|x(t)+y(t)|\rightarrow c, \, t\rightarrow\infty.$$
If $c=0,$ then the anti-synchronization is called complete
anti-synchronization. Anti-synchronization can also be bone using
drive and response systems of the same type, or different types of
chaotic systems of fractional differential equations, as above. We
will con\-sider both theses types of anti-synchronization using
the systems (\ref{T1}) and (\ref{R1}) \cite{Seb}.

The anti-synchronization between two identical fractional T
systems is done in a simi\-lar manner like in the synchronization
case. We consider the drive system (\ref{T1}) and the response
system with a control $u(t)=[u_1(t),u_2(t),u_3(t)]^T$ (\ref{T2}).

In this case the error system is given by adding  (\ref{T2}) and
(\ref{T1}), and we get
\begin{equation}\label{Te11}
             \left\{%
             \begin{array}{ll}
                D^{\alpha_1}_t e_1(t)=a_1(e_2(t)-e_1(t))+u_1(t),\\
                D^{\alpha_2}_t e_2(t)=(c_1-a_1)e_1(t)-a_1(x_2(t)z_2(t)+x_1(t)z_1(t))+u_2(t),\\
                D^{\alpha_3}_t e_3(t)=-b_1e_3(t)+x_2(t)y_2(t)+x_1(t)y_1(t)+u_3(t),\\
                \end{array}
                \right.
\end{equation}where $e_1(t):=x_2(t)+x_1(t),$
$e_2(t):=y_2(t)+y_1(t),$ $e_3(t):=z_2(t)+z_1(t).$ Now we will
re-write the control $u$ by choosing a suitable control  $v$ as a
function depending on the error states $e_1, e_2, e_3.$ We
re-refine $u$ as
\begin{equation}\label{c11}
\begin{array}{ll}
u_1(t):=v_1(t),\\
u_2(t):=a_1(x_2(t)z_2(t)+x_1(t)z_1(t))+v_2(t),\\
u_3(t):=-x_1(t)y_1(t)-x_2(t)y_2(t)+v_3(t).
\end{array}
\end{equation}
Substituting (\ref{c11}) in (\ref{Te11}) we obtain
\begin{equation}\label{Te21}
             \left\{%
             \begin{array}{ll}
                D^{\alpha_1}_t e_1(t)=a_1(e_2(t)-e_1(t))+v_1(t),\\
                D^{\alpha_2}_t e_2(t)=(c_1-a_1)e_1(t)+v_2(t),\\
                D^{\alpha_3}_t e_3(t)=-b_1e_3(t)+v_3(t).\\
                \end{array}
                \right.
\end{equation}
Two considered systems are globally synchronized if relation the
error to converge to 0, when $t\rightarrow\infty.$ We have to
choose the feedback control $v$ such that relation (\ref{v}) is
fulfilled. One choice of $A$ is $$A:=\begin{bmatrix}
  0 & a_1 & 0\\
  -(c_1-a_1) & c_1 & 0\\
  0 & 0 & 2b_1 \\
\end{bmatrix}$$ and the feedback functions are
\begin{equation}\label{vT1}
\begin{array}{ll}
v_1(t)=a_1e_2(t),\\
v_2(t)=-(c_1-a_1)e_1(t)+c_1e_2(t),\\
v_3(t)=2b_1e_3(t).
\end{array}
\end{equation}
The error system (\ref{Te21}) becomes
\begin{equation}\label{Te31}
             \left\{%
             \begin{array}{ll}
                D^{\alpha_1}_t e_1(t)=a_1e_1(t),\\
                D^{\alpha_2}_t e_2(t)=c_1e_2(t),\\
                D^{\alpha_3}_t e_3(t)=b_1e_3(t).\\
                \end{array}
                \right.
\end{equation}

\begin{proposition}
The systems of fractional differential equations \emph{(\ref{T1})}
and \emph{(\ref{T2})} will a\-pproach global asymptotical
anti-synchronization for any initial conditions and with the
feedback control \emph{(\ref{vT1})}.
\end{proposition}
\textbf{Proof:} The proof is similar to that of Proposition
\ref{Pr1}. \hfill $\Box$\\

In case of  R\"{o}ssler system as drive system (\ref{R1}) and T
system (\ref{T1}) as response system, the error system is obtained
by adding (\ref{T11}) and (\ref{R1}) and it has the form
\begin{equation}\label{Re11}
             \left\{%
             \begin{array}{ll}
                D^{\alpha_1}_t e_1(t)=-a_1e_1(t)+(a_1-1)e_2(t)-e_3(t)+a_1(x_1(t)-y_1(t))
                +y_2(t)+z_2(t)+u_1(t),\\
                D^{\alpha_2}_t
                e_2(t)=(c_1-a_1+1)e_1(t)+a_2e_2(t)-(c_1-a_1)x_1(t)-x_2(t)-
                a_2x_2(t)z_2(t)- \\ \quad \quad \quad \quad \quad \,  -a_2y_2(t)+u_2(t),\\
                D^{\alpha_3}_t
                e_3(t)=-(b_1+c_2)e_3(t)+b_2+x_2(t)y_2(t)+b_1z_1(t)+z_1(t)x_1(t)+
                c_2z_2(t)+u_3(t),\\
                \end{array}
                \right.
\end{equation}where $e_1(t):=x_2(t)+x_1(t),$
$e_2(t):=y_2(t)+y_1(t),$ $e_3(t):=z_2(t)+z_1(t).$ We re-write
control $u$ by considering the control $v$ such that
\begin{equation}
\begin{array}{ll}
u_1(t):=-a_1(x_1(t)-y_1(t))-y_2(t)-z_2(t)+v_1(t),\\
u_2(t):=(c_1-a_1)x_1(t)+x_2(t)+a_2x_2(t)z_2(t)+a_2y_2(t)+v_2(t),\\
u_3(t):=-b_2-x_2(t)y_2(t)-b_1z_1(t)-z_1(t)x_1(t)-c_2z_2(t)+v_3(t),
\end{array}
\end{equation}and the error system (\ref{Re11}) becomes
\begin{equation}\label{Re21}
             \left\{%
             \begin{array}{ll}
                D^{\alpha_1}_t e_1(t)=-a_1e_1(t)+(a_1-1)e_2(t)-e_3(t)+v_1(t),\\
                D^{\alpha_2}_t e_2(t)=(c_1-a_1+1)e_1(t)+a_2e_2(t)+v_2(t),\\
                D^{\alpha_3}_t e_3(t)=-(b_1+c_2)e_3(t)+v_3(t).\\
                \end{array}
                \right.
\end{equation}
As in the case above, we have to chose the feedback control $v$
such that the relation (\ref{v}) to be fulfilled. One choice of
matrix $A$ will be
$$A:=\begin{bmatrix}
  2a_1 & -(a_1-1) & 1 \\
  -(c_1-a_1+1) & 0 & 0\\
  0 & 0 & 2b_1+c_2 \\
\end{bmatrix}$$ and the feedback functions are
\begin{equation}\label{vR1}
\begin{array}{ll}
v_1(t)=2a_1e_1(t)-(a_1-1)e_2+e_3(t),\\
v_2(t)=-(c_1-a_1+1)e_1(t),\\
v_3(t)=(2b_1+c_2)e_3(t).
\end{array}
\end{equation}
The error system (\ref{Re21}) becomes
\begin{equation}\label{Re31}
             \left\{%
             \begin{array}{ll}
                D^{\alpha_1}_t e_1(t)=a_1e_1(t),\\
                D^{\alpha_2}_t e_2(t)=a_2e_2(t),\\
                D^{\alpha_3}_t e_3(t)=b_1e_3(t).\\
                \end{array}
                \right.
\end{equation}

\begin{proposition}
The systems of fractional differential equations \emph{(\ref{R1})}
and \emph{(\ref{T11})} will a\-pproach global asymptotical
synchronization for any initial conditions and with the feedback
control \emph{(\ref{vR1})}.
\end{proposition}
\textbf{Proof:} The proof is similar to that of Proposition
\ref{Pr1}. \hfill $\Box$

\section{Secure information using synchronized chaotic systems of fractional
differential equations}

In the following we will do the encryption and decryption of a
message using coupled chaotic systems of fractional differential
equations and their synchronization. For this we use two coupled
systems, for example (\ref{T1}) as drive system and (\ref{T2}) as
response. They have a chaotic behaviour and are also synchronized
and anti-synchronized. The sender uses the fractional system
(\ref{T1}) and the receiver uses (\ref{T2}). They both choose
values for variables $z_1(t)$ and $z_2(t)$ as public keys, after a
period of time, after the synchronization between the considered
fractional systems took place \cite{Ch}, \cite{GH}. The fractional
derivative is also sent as a public key to the receiver.
\par The message that one part wants to send is called
\emph{plaintext} and its correspondent by decryption is called
\emph{ciphertext}. The plaintext and the ciphertext are
represented by numbers, each letter from the alphabet is replaced
by a corresponding number. So instead of letters and numbers we
will use numbers from 0 to 35, in this case. The most general case
is to consider small letters, capital letters and special
characters. To each of it, a corresponding number is associated,
in ASCII representation. In this case the formula corresponding
for encryption, (respectively for decryption) is given by:
\begin{equation}
\begin{array}{ll}
 c_i :=p_i + k_i\, mod(36), \\
 p_i := c_i - k_i\, mod(36),\\
\end{array}
\end{equation}
where $k_i$ are public keys that mask the message. For each letter
or number we use a randomly generated key $\{k_1, k_2, . . .
k_n\}.$ Actually, each key $k_j$ hides the piece of message $p_j.$

We consider the example message  "Hello Oscar". This is
represented in the following table, each piece of message with the
corresponding randomly generated key, and in the next table the
plaintex, the cipertext and the message received:
\\ \\
\begin{tabular}{|c|c|c|c|c|c|c|c|c|c|c|c|c|c|c|c|c|}
  \hline
 Key& $k_1$ & $k_2$ & $k_3$ & $k_4$ & $k_5$ & $k_6$ & $k_7$ & $k_8$ & $k_9$ & $k_{10}$ \\
  \hline
 Message & h & e & l & l & o & o & s & c & a & r \\
  \hline
\end{tabular}\\
\\
\par Beginning with $t_0$, in our example equal with 1300, both dynamical systems are
synchronized, so $t$ takes values greater than $t_0=1300,$ in the
synchronized state. The first system sends the encoded message,
letter by letter, and for each letter a key is randomly generated.
The second system receives the encrypted message and also the key
with which is was encrypted. In this example we work with the
third component of the considered fractional systems, $z_1(t)$ and
with the keys  $k_i, \, i=\overline{1, 10}.$ Decryption is done
using the third component of the second system $z_2(t),$ with
$t\geq t_0,$ as it is illustrated in Table 2.
\\ \\
\begin{table}\label{table1}
\caption{Ciphertext}
\begin{tabular}{|c|c|c|c|c|c|c|c|c|c|c|c|c|c|c|c|c|}
  \hline
 Time $t$& $z_1(t)$ & Key $k$ & Plaintext $p$ & Ciphertext  $c=p+k\, mod(36)$ \\
  \hline
  1301 & 4416 &  422711... & h(18) & 0 \\
  1302 & 4433 &  312595... & e(15) & 33 \\
  1303 & 4449 &  155104... & l(22) & 15 \\
  1304 & 4466 &  974929... & l(22) & 15 \\
  1305 & 4483 &  271362... & o(25) & 9 \\
  1306 & 4500 &  333989... & o(25) & 9 \\
  1307 & 4518 &  649588... & s(29) & 14 \\
  1308 & 4535 &  799296... & c(13) & 19 \\
  1309 & 4552 &  379435... & a(11) & 17 \\
  1310 & 4569 &  963282... & r(28) & 22 \\
  \hline
\end{tabular}
\end{table}
\begin{table}
\caption{Plaintext}
\begin{tabular}{|c|c|c|c|c|c|c|c|c|c|c|c|c|c|c|c|c|}
  \hline
 Time $t$& $z_2(t)$ & Key $k$  & Ciphertext $p$ & Plaintext $p=c-k mod(36)$\\
   \hline
  1301 & 4416 &  422711... & 0 & 18(h) \\
  1302 & 4433 &  312595... & 33 & 15(e) \\
  1303 & 4449 &  155104... & 15 & 22(l) \\
  1304 & 4466 &  974929... & 15 & 22(l) \\
  1305 & 4483 &  271362... & 9 & 25(o) \\
  1306 & 4500 &  333989... & 9 & 25(o) \\
  1307 & 4518 &  649588... & 14 & 29(s) \\
  1308 & 4535 &  799296... & 19 & 13(c) \\
  1309 & 4552 &  379435... & 17 & 11(a) \\
  1310 & 4569 &  963282... & 22 & 28(r) \\
  \hline
\end{tabular}
\end{table}

\section{Conclusions}

In this paper we have presented chaotic behavior of some
fractional systems and we have presented a way of encryption and
decryption a message. These techniques can be done using other
chaotic systems of fractional differential equations, like Lorenz,
Chua, hyperchaotic R\"{o}ssler systems and pair of these
fractional systems, but also on some systems of fractional
differential equations endowed with another structure, such as
metriplectic structure, (almost) Leibniz structure \cite{Ch2}. A
suitable control can be chosen to achieve synchronization or
anti-synchronization.

For a suitable control $u,$ respectively $v,$ synchronization and
anti-synchronization behave in the same manner, as the error
system of fractional differential equations coincide.

Encryption and decryption was done here letter by letter, using a
Maple 11 program. In a discrete case the procedure is more
"commercial" because encoding and decoding the message is done for
the entire message at the same time.

Another approach for synchronization has been studied and
presented in my PhD thesis \cite{Ch1}.

\end{document}